\documentclass[12pt,twoside]{amsart}

  \usepackage{url}
  \usepackage{amssymb,amsmath}
     \usepackage{graphicx}

   \usepackage[a4paper, margin=2cm]{geometry}

  \input xypic
  \xyoption{all}

 \usepackage{tikz}
\usetikzlibrary{arrows}

 \usepackage{times}

  \usepackage{amsmath, amsthm, amsfonts}

\def \Sinverse  { U_{I,\epsilon}(M^{+}) }

\newtheorem{theorem}{Theorem}[section]

\newtheorem{corollary}[theorem]{Corollary}
\newtheorem{definition}[theorem]{Definition}

\newtheorem{lemma}[theorem]{Lemma}

\def \fa {\mathfrak{a}}
\def \p{\partial}

\def \bR {\mathbb R}
\def \bh {\mathbb H}

 \def \bk {\bf k}

               \def \bR{{\mathbb R}}

\def \p{\partial}

\def \v {\vskip 0.1in}
\def \n {\noindent}

\begin{document}

\begin{center}
  {\LARGE Contact Invariants, Open String Invariants\\ and Weinstein Conjecture}
  \end{center}

  \noindent
  \begin{center}
   {\large An-Min Li}\footnote{partially supported by a NSFC grant}\\[5pt]
      Department of Mathematics, Sichuan University\\
        Chengdu, PRC\\[5pt]
  {\large  Li Sheng}\footnote{partially supported by a NSFC grant}\\[5pt]
      Department of Mathematics, Sichuan University\\
        Chengdu, PRC\\[5pt]
\end{center}

\begin{abstract}
In this note we propose a theory of contact invariants and open string invariants, which can be used to study the Weinstein Conjecture.
\end{abstract}

 % \tableofcontents
\v\v

In \cite{LS} we propose a theory of contact invariants and open string invariants, assuming that every periodic orbit of the Reeb vector field is either non-degenerate or of Bott-type, where we choose   the complex structure $\tilde{J}$ such that $L_X\tilde{J}=0$ on periodic orbits. In this note we do not choose the complex structure $\tilde{J}$ such that $L_X\tilde{J}=0$ on periodic orbits, so we don't mod the $S^1$ action on every periodic orbit. Thus we have to consider the moduli space of partially decorated stable nodal surfaces ( see \cite{BEHWZ}). Then the moduli space of perturbed $J$ holomorphic maps has codimension 1 boundary. But we can still define the contact invariants and open string invariants, because the integrals  \eqref{c_invariant}, \eqref{l_invariant} are convergent. Thanks of the exponential decay estimates of the gluing maps with respect to the gluing parameter.
\v

\section{\bf Symplectic manifolds with cylindrical ends}
\v

\subsection{Contact manifolds}

Let $(\widetilde{M},\lambda)$ be a $(2n-1)$-dimensional compact manifold
equipped with a contact form $\lambda$. We recall that a contact form
$\lambda$ is a 1-form on $\widetilde{M}$ such that $\lambda\bigwedge
(d\lambda)^{n-1}$ is a volume form. Associated to $(\widetilde{M},\lambda)$
we have the contact structure $\xi=\ker(\lambda)$, which is a
$(2n-2)$-dimensional subbundle of $T\widetilde{M}$, and $(\xi, d\lambda|_{\xi})$
defines a symplectic vector bundle. Furthermore, there is a unique
nonvanishing vector field $X=X_{\lambda}$, called the Reeb vector field,
defined by the condition
$$\lambda(X)=1,\;\;\;\;i_Xd\lambda=0.$$
We have a canonical splitting of $T\widetilde{M}$,
$$ T\widetilde{M}=\mathbb R X\oplus\xi,$$
where $\mathbb R X$ is the line bundle generated by $X$.

\v

\subsection{Cylindrical almost complex structures}\label{subsection:7.2}

Let
\begin{equation}\label{cylinder}
M^+=M^{+}_{0}\bigcup\left\{[0,\infty)\times \widetilde{M}\right\}\end{equation}
be a manifold with cylindrical end,
where $M^{+}_{0}$ is a compact manifold with boundary, $\widetilde{M}$ is a compact contact manifold with contact form $\lambda$.
Let
$${\bf \Phi}=\left\{\phi \in C^{\infty}(\mathbb R, [0,1]) | \phi'\geq 0 \right\}$$
and define the 1-form $\lambda_{\phi}=\phi(a)\lambda$ over $[0,\infty)\times \widetilde{M}$.
Suppose that $M_0^{+}$ has a symplectic form $\omega$. Denote by $\omega_{\phi}$  the symplectic form of $M^{+}$ such that
$\omega_{\phi}|_{M_0^{+}}=\omega$, and over the cylinder $[0,\infty)\times \widetilde{M}$
\begin{equation}\omega_{\phi} =
\phi d\lambda + \phi^{\prime}d a\wedge \lambda.
\end{equation} We will also consider $\mathbb R\times \widetilde{M}$. Denote by $N$ one of $M^+$ and $\mathbb R\times \widetilde{M}$.
\v
We choose a $d \lambda$ compatible almost complex structure $\widetilde{J}$ for the symplectic vector bundle
$(\xi, d\lambda)\rightarrow \widetilde{M}$ such that
\begin{equation}\label{metric}
g_{\widetilde{J}(x)}(h,k)=d\lambda(x)(h,\widetilde{J}(x)k)
\end{equation}
for all $x\in \widetilde{M},\;h,k\in\xi_x$, defines a smooth fibrewise
metric for $\xi$. Denote by $\Pi:T\widetilde{M}\rightarrow\xi$ the projection along $X$.
We define a Riemannian metric $\langle\;,\;\rangle$ on $\widetilde{M}$ by
\begin{equation}\label{metric1}
\langle h,k\rangle=\lambda(h)\lambda(k)+g_{\widetilde{J}}(\Pi h,\Pi k)
\end{equation}
for all $h,k\in T\widetilde{M}$.
\vskip 0.1in
Given a $\widetilde{J}$ as above there is an associated almost complex
structure $J$ on $\mathbb R\times \widetilde{M}$ defined by
\begin{equation}\label{complex structure}
J\mid_{\xi}= \tilde{J},\;\;\;JX=-\frac{\partial}{\partial a},\;\;J(\frac{\partial}{\partial a})=X,
\end{equation}
where $a$ is the canonical coordinate in $\mathbb R$.

It is easy to check that $J$ defined by \eqref{complex structure}
is  $\omega_{\phi}$-tame over the cylinder end. We can choose an almost complex structure $J$ on $M^+$ such that $J$ is tamed by $\omega$ and over the cylinder end  $J$ is given by \eqref{complex structure}.
\v
 \noindent For any $\phi \in \Phi $
\begin{equation}\label{omega_forms}
\langle v,w\rangle_{\omega_{\phi}} = \frac{1}{2}\left(
\omega_{\phi} (v,Jw) + \omega_{\phi} (w,Jv)\right) \;\;\;\;\;
\forall \;\; v, w \in TN
\end{equation}
defines a Riemannian metric
on $N$.  Note that $\langle \;,\;\rangle_{\omega_{\phi}}$ is not
complete. We choose another metric $( \;,\;)$ on $N$
such that
\begin{equation}\label{omega_forms_on_M0}( \;,\;) = \langle
\;,\;\rangle_{\omega_{\phi}} \;\;\;\;on \;\; M^{+}_{0}
\end{equation}
 and over
the tubes
\begin{equation} \label{omega_forms_on_tubes}
((a,v),(b,w))= ab + \lambda (v)\lambda
(w) + g_{\widetilde{J}}(\Pi v, \Pi w).
\end{equation}
It is easy to see that $( \;,\;)$ is a complete
metric on $N$. \vskip 0.1in \noindent
\v\v
\subsection{Neighbourhoods of Lagrangian submanifolds}\label{neighbor_L}

Let $(M,\omega)$ be a compact symplectic manifold, $L\subset M$ be a compact Lagrangian submanifold.
The following Theorem is well-known.
\begin{theorem}\label{Lagrangian neighbourhood theorem}
 Let $(M,\omega)$ be a symplectic manifold of dimension
$2n$, and $L$ be a compact Lagrangian submanifold. Then there exists a neighbourhood $U\subset T^*L$ of the zero section, a neighbourhood $V \subset M$ of $L$, and a diffeomorphism $\phi:U\rightarrow V$ such that
\begin{equation}
\phi^{*}\omega=-d \Lambda,\;\;\;\; \phi|_{L}=id,
\end{equation}
where $\Lambda$ is the canonical Liouville  form.
\end{theorem}

 Let $(x_1,\cdots, x_n)$ be a
local coordinate system on $O\subset L$, there is a canonical coordinates
$$(x_1,\cdots,x_n,y_1,\cdots,y_n)$$ on $T^{*}O=T^{*}L|_{O}.$ In terms of this coordinates
the  Liouville form can be written as
$$\Lambda =\sum y_{i}dx_{i}.$$
Let $\pi :T^{*}L \rightarrow L$ be the canonical projection.
Suppose that given a Riemannian metric on $L$, in terms of the coordinates $x_1,...,x_n$, $g_{L}=\sum\limits_{i,j=1}^{n} g_{ij}dx_idx_j.$ It naturally induced a metric on $T^{*}L.$

\v

Denote by $S^{n-1}(1)$ (resp.$B_{1}(0)$) the Euclidean unit sphere (resp. the Euclidean unit ball). Consider the coordinates transformation between the sphere coordinates and the Cartesian coordinate
\begin{align}
\Psi:(0,1]\times S^{n-1}(1)&\rightarrow  B_{1}(0) \nonumber\\
(r,\theta_{1},\cdots,\theta_{n-1})&\rightarrow   (y_1,\cdots,y_{n}).
\end{align}
Consider the unit sphere bundle $\widetilde{M}$ and the unit ball bundle $\mathbb D_{1}(T^{*}L)$ in $T^{*}L,$ in terms of the coordinates $(x_{1},\cdots,x_{n},y_{1},\cdots,y_n)$
\begin{align}
 \widetilde{M}|_{\pi^{-1}(O)}&=\{ (x_{1},\cdots,x_{n},y_{1},\cdots,y_n)\in \pi^{-1}(O)\;|\; \sum_{i,j=1}^{n} g^{ij}(x)y_{i}y_{j}=1 \},\\
\mathbb D_{1}(T^{*}L)|_{\pi^{-1}(O)}&=\{ (x_{1},\cdots,x_{n},y_{1},\cdots,y_n)\in \pi^{-1}(O)\;|\; \sum_{i,j=1}^{n} g^{ij}(x)y_{i}y_{j}\leq 1 \}.
\end{align}
Denote $\lambda= -\Lambda \mid_{\widetilde{M}}$.
We have
$$\Lambda=-|{\bf y}|\lambda,$$ where $|\cdot|$ denotes the Euclidean norm. $\lambda $ is a contact form, i.e., $(\widetilde{M}, \lambda)$ is a contact manifold. Put $\xi=\ker(\lambda)$. Denote
$X=- \sum g^{ij}y_i\frac{\partial}{\partial x_j}\mid_{\widetilde{M}}$. Then $X $ is the Reeb vector field. Denote $z=1-|{\bf y}|.$ By Theorem \ref{Lagrangian neighbourhood theorem} we consider $M-L$ as
$$M^+= M_0^{+}\bigcup\{(0,1]\times \widetilde{M}\}$$
with the symplectic form
\begin{equation}\omega =- d\Lambda=
 (1-z) d\lambda - d z \wedge\lambda,
\end{equation}
where $M^+:= M-L$ and $M^{+}_{0}$ is a compact symplectic manifold with boundary.
\v
We choose the  {\em neck stretching  technique}. Let $\phi :[0,\infty)\rightarrow [0,1)$ be a smooth
function satisfying, for any $k>0,$
$$\phi^{\prime}>0,
\;\;\lim_{a\to \infty}\phi(a)= 1, \;\;\phi(0)= 0. $$
Through $\phi$ we consider
$M^+$ to be $M^{+} = M_{0}^{+}\bigcup\{[0,\infty)\times
\widetilde{M}\} $ with symplectic form
$\omega_{\phi}|_{M_0^{+}}=\omega$, and over the cylinder $[0,\infty)\times \widetilde{M}$
\begin{equation}\omega_{\phi} =- d\Lambda=
(1-\phi) d\lambda - \phi^{\prime}d a\wedge \lambda.
\end{equation}
Denote
$${\bf \Phi}=\left \{ \phi :[0,
\infty )\rightarrow [0, 1) | \phi^{\prime}
>0 \right \}.$$

We choose a $d \lambda$ compatible almost complex structure $\widetilde{J}$ for the symplectic vector bundle
$(\xi, d\lambda)\rightarrow \widetilde{M}$.
There is an associated almost complex
structure $J$ on $\mathbb R\times \widetilde{M}$ defined by
\begin{equation}\label{complex structure}
J\mid_{\xi}= \tilde{J},\;\;\;JX= \frac{\partial}{\partial a},\;\;J(\frac{\partial}{\partial a})=-X,
\end{equation}
where $a$ is the canonical coordinate in $[0,\infty)$.

\section{\bf $J$-holomorphic maps with finite energy}
\v

Let $(\Sigma,i)$ be a
compact Riemann surface and $P\subset\Sigma$ be a finite
collection of puncture points. Denote $\stackrel{\circ}{\Sigma}
=\Sigma\backslash P.$ Let ${u}:\stackrel{\circ}{\Sigma}
\rightarrow N$ be a ${J}$-holomorphic map, i.e., ${u}$ satisfies
\begin{equation}\label{j_holomorphic_maps}
d{u}\circ i={J}\circ d{u}.
\end{equation}
We write $u=(a, \widetilde{u})$ and define
\begin{equation}\label{definition_of_energy_on_complex_manifolds}
\widetilde{E}(u)=\int_{{\Sigma}}\widetilde{u}^{\ast}
 d\lambda .
\end{equation}
For any $J$-holomorphic map
$u:\stackrel{\circ}{\Sigma}\rightarrow N$ and any $\phi \in \Phi $
the energy $E_{\phi}(u)$ is defined by
\begin{equation}\label{definition_of_energy}
E_{\phi}(u)=\int_{ {\Sigma}}u^{\ast}\omega_{\phi}.
\end{equation}
Let $z=e^{s+2 \pi \sqrt{-1}t}.$ One computes
\begin{equation}\label{omega_forms_cylinder}
u^{\ast}\omega_{\phi}=
(\phi d\lambda \left((\pi\widetilde{u})_s,
(\pi\widetilde{u})_t)\right) + {\phi}^{\prime}(a^{2}_s + a^{2}_t
))ds\wedge dt,
\end{equation} which is a nonnegative integrand.
Following \cite{HWZ1} we
impose an energy condition on $u$.
A $J$-holomorphic map $u:\stackrel{\circ}{\Sigma} \rightarrow N $
is called a finite energy $J$-holomorphic map if over the cylinder end
\begin{equation}\label{finite_energy_j_holomorphic_maps}
\sup_{\phi \in
\Phi }\left \{\int_{{\Sigma}}u^{*} \omega_{\phi}
\right \}+\int_{\Sigma}u^{*}d\lambda <\infty.
\end{equation}
  For a $J$-holomorphic
map $u:{\Sigma} \rightarrow {\mathbb
R}\times\widetilde{M}$  we write $u=(a, \widetilde{u})$ and define
\begin{equation}\label{definition_of_energy_on_complex_manifolds}
\widetilde{E}(u)=\int_{{\Sigma}}\widetilde{u}^{\ast}
 d\lambda .
\end{equation}
Denote
$$\widetilde{E}(s)=\int_s^{\infty}\int_{S^1}\widetilde{u}^{\ast}(
d\lambda).$$ Then $$\widetilde{E}(s)=\int_s^{\infty}\int_{S^1} |\Pi
\widetilde{u}_t |^2dsdt,$$
\begin{equation}\label{deriative_of_energy_on_Z}
\frac{d\widetilde{E}(s)}{ds}=-\int_{S^1} |\Pi\widetilde{u}_t |^2dt.
\end{equation}
%where $\|\cdot\|$ denotes the norm defined by \eqref{omega_forms_on_tubes}.
Here and later we use $ |\cdot |$ denotes the norm with respect to the metric defined by \eqref{omega_forms_on_tubes}.

 \v\n
Following Hofer et al. \cite{BEHWZ}  we assume that
\v
{\bf Condition A.} the almost complex structure $J$   either nondegenerate or of Bott-type.
\v
The following theorems are well-known:
\v
\begin{theorem} \label{main theorem}
Denote $\mathbb D_1(0)=\{z\in \mathbb C|\; |z|< 1\}$. Assume that {\bf Condition A} holds.  Let $u=(a,\widetilde{u}):\mathbb D_1(0) \rightarrow {\mathbb R} \times
\widetilde{M}$ be a nonconstant $J$-holomorphic map with finite
energy. Put $z=e^{-s+2\pi \sqrt{-1}t}$. Then
$$\lim_{s\rightarrow \infty}\widetilde
u(s,t)=x(kTt) $$ in $C^{\infty}(S^1)$ for some $kT$-periodic orbit
$x$ of the Reeb vector field.
\end{theorem}
Following Hofer (see \cite{HWZ1}) we introduce a convenient local coordinates $(\vartheta, w^1,...,w^{2n-2})$ around the periodic orbit $x$, we call it a pseudo-Darboux coordinate system,
such that
\begin{equation}
\lambda = f\lambda_0,\end{equation}
where $\lambda_{o}=d\vartheta+\sum_{i=1}^{n-1} w^{i}dw^{n-1+i}$ and $f:U\rightarrow \mathbb R$ is a smooth function satisfying
\begin{equation}\label{standard_lambda-1}
f(\vartheta, 0)=T,\;\;df(\vartheta,0)=0
\end{equation}
for all $\vartheta\in S^1$.
\v
\begin{theorem}\label{partial_a_theta_convergence_zero-1} Assume that {\bf Condition A} holds.  Let $u=(a,\widetilde{u}):\mathbb D_1(0) \rightarrow {\mathbb R}  \times \widetilde{M}$ be as in Theorem \ref{main theorem}.
Then there are constants $\ell_0$,
$\vartheta_0$ and $0<\mathfrak{c}<\frac{1}{2}$  such that for all ${\bf n}=(n_1,n_2)\in {\mathbb Z_{\geq 0}^2 }$
\begin{eqnarray}
\label{exponential_decay_a}
|\partial^{\bf n}[a(s,t)-kTs-\ell_0]|\leq \mathcal C_{\bf n} e^{-\mathfrak{c}|s|}\\
\label{exponential_decay_theta}
|\partial^{\bf n}[\vartheta(s,t)-kt-\vartheta_0]|\leq \mathcal C_{\bf n} e^{-\mathfrak{c}|s|} \\
\label{exponential_decay_y}|\partial^{\bf n} {\mathbf w}(s,t)|\leq \mathcal C_{\bf n}
e^{-\mathfrak{c} |s|},
\end{eqnarray}
 where $\mathcal{C}_{\bf n}$ are constants.
\end{theorem}
\v\v

\section{\bf Weighted sobolev norms}\label{weight_norm}

\v
Consider ${\mathbb{R}}\times \widetilde{M}$ and $M^+=M^{+}_{0}\bigcup\left\{[0,\infty)\times \widetilde{M}\right\}.$ Let $N$ be one of ${\mathbb{R}}\times \widetilde{M}$ and $M^+$. Suppose that  $\Sigma= \bigcup\limits_v \Sigma_{v}$ is Riemann surface with nodal points $\{q_{1},\cdots,q_{\mathfrak{I}}\}$, puncture points $\{p_{1},\cdots,p_{\nu}\}$ and  $u:{\Sigma}\rightarrow \bigcup N_i$ is a continuous map such that the restriction of $u$ to each smooth component is smooth, where $\bigcup N_i$ denotes the union of some copy of $N$. We choose cylinder coordinates $(s,t)$ on $\Sigma$ near each nodal point and each puncture point. We choose a local pseudo-Darboux coordinate system near each periodic orbit on $N$. Let $\stackrel{\circ}{\Sigma}=\Sigma-\{q_{1},\cdots,q_{\mathfrak{I}},p_{1},\cdots,p_{\nu} \}$.
\v

Over each tube the linearized operator
$D_{u}$ takes the following form
\begin{equation}
D_{u}=\frac{\partial}{\partial
s}+J_0\frac{\partial} {\partial t}+S = \bar{\partial}_{J} +
S.
\end{equation}
By exponential decay we have
$$
\left|\frac{\partial^{k+l} }{\partial s^k \partial t^{l} }S\right|\leq C_{k,l}e^{- \mathfrak{c} s}
$$
for some constant $C_{k,l}>0$ for $s$ big enough.
Therefore, the operator $H_s=J_0\frac{d}{dt}+S$ converges
to $H_{\infty}=J_0\frac{d}{dt}$. Obviously, the operator $D_u$ is
not Fredholm operator because over each   puncture and node the
operator $H_{\infty}=J\frac{D}{dt}$ has zero eigenvalue. The $\ker
H_{\infty}$ consists of constant vectors. To recover a Fredholm theory we use
weighted function spaces. We choose a weight $\alpha$ for each
end. Fix a positive function $W$ on $\Sigma$ which has order equal
to $e^{\alpha |s|} $ on each end, where $\alpha$ is a small
constant such that $0<\alpha<\mathfrak{c} $ and over each end
$H_{\infty}- \alpha = J_0\frac{d}{dt}- \alpha $ is invertible. We
will write the weight function simply as $e^{\alpha |s|}.$ Denote by $C(\Sigma;u^{\ast}T(\bigcup N_i))$ all tangent vector fields $h$ on $\bigcup N_i$ along $u$ satisfying
\v
{\bf (a)} $h\in C^{0}({\Sigma},u^*T(\bigcup N_i))$,
\v
{\bf (b)} the restriction of $h$ to each smooth component is smooth.
\v
 For any
section $h\in C(\Sigma;u^{\ast}T(\bigcup N_i))$ and section $\eta \in
\Omega^{0,1}(u^{\ast}T(\bigcup N_i))$ we define the norms
\begin{eqnarray}\label{def_1_p_alpha}
&&\|h\|_{1,p,\alpha}=\sum_{v}\left(\int_{\Sigma_v}( |h|^p+ |\nabla
h|^p)d\mu\right)^{1/p}
+\sum_{v}\left(\int_{\Sigma_v}e^{2\alpha|s|}(|h|^2+|\nabla h|^2)
d\mu\right)^{1/2} \\ \label{def_p_alpha}
&&\|\eta\|_{p,\alpha}=\sum_{v}\left(\int_{\Sigma_v}|\eta|^p d\mu
\right)^{1/p}+ \sum_v\left(\int_{\Sigma_v}e^{2\alpha|s|}|\eta|^2
d\mu\right)^{1/2}
\end{eqnarray} for $p\geq 2$, where all norms and
covariant derivatives are taken with respect to the  metric
$(\;,\;)$ on $u^{\ast}T(\bigcup N_i)$ defined in \eqref{omega_forms_on_tubes}, and
the metric on $\stackrel{\circ}{\Sigma}$. Denote
\begin{eqnarray}
&&{\mathcal C}(\Sigma;u^{\ast}T(\bigcup N_i))=\{h
\in C(\Sigma;u^{\ast}T(\bigcup N_i)); \|h\|_{1,p,\alpha}< \infty
\},\\
&&{\mathcal C}(u^{\ast}T(\bigcup N_i)\otimes \wedge^{0,1})
=\{\eta\in \Omega^{0,1}(u^{\ast}T(\bigcup N_i)); \|\eta\|_{p,\alpha}< \infty
\}.
\end{eqnarray}
 Denote by $W^{1,p,\alpha}(\Sigma;u^{\ast}T(\bigcup N_i))$ and
$L^{p,\alpha}(u^{\ast}T(\bigcup N_i)\otimes \wedge^{0,1})$ the completions of
${\mathcal C}(\Sigma;u^{\ast}T(\bigcup N_i))$ and ${\mathcal C}(u^{\ast}T(\bigcup N_i)\otimes
\wedge^{0,1}) $ with respect to the norms \eqref{def_1_p_alpha} and \eqref{def_p_alpha}
respectively. Then the operator $D_u:   W^{1,p,\alpha}\rightarrow L^{p,\alpha}$
is a Fredholm operator.
\vskip 0.1in
\noindent
\v
For each bounded nodal $q_{i},$ denote  $\bh_{q_{i}}=T_{q_{i}}N$, for each unbounded nodal $q_{i},$ denote  $\bh_{q_{i}}=T_{q_{i}}\widetilde{M}\oplus (span\{\frac{\p}{\p a}\}.$
Put
$$\bh= \left(\oplus_{j=1}^\nu (T_{p_{j}}\widetilde{M}\oplus (span\{\frac{\p}{\p a}\})\right) \bigoplus \left(\oplus_{i=1}^{\mathfrak{I}} \bh_{q_{i}}\right) ),$$
 $$h_{0}=(h_{1 0},...,h_{\nu 0},h_{(1+\nu)0},...,h_{  (\mathfrak{I}+\nu) 0}).$$
$h_0$ may be considered as a vector field in the coordinate neighborhood.
We fix a cutoff function $\rho$:
\[
\rho(s)=\left\{
\begin{array}{ll}
1, & if\ |s|\geq d, \\
0, & if\ |s|\leq \frac{d}{2}
\end{array}
\right.
\]
where $d$ is a large positive number. Put
$$\hat{h}_0=\rho h_0.$$
Then for $d$ big enough $\hat{h}_0$ is a section in $C^{\infty}(\Sigma; u^{\ast}TN)$
supported in the tube $\{(s,t)||s|\geq \frac{d}{2}, t \in S ^1\}$.
Denote
$${\mathcal W}^{1,p,\alpha}=\{h+\hat{h}_0 | h \in
W^{1,p,\alpha},h_0 \in \bh\}.$$
\vskip 0.1in
\noindent
 We define the weighted Sobolev  norm  on ${\mathcal W}^{1,p,\alpha}$ by $$\|( h, \hat{h}_{0})\|_{\Sigma,1,p,\alpha}=
 \|h\|_{\Sigma,1,p,\alpha} + |h_{0}| .$$
Obviously, the operator $D_u:   \mathcal W^{1,p,\alpha}\rightarrow L^{p,\alpha}$
is also a Fredholm operator.
\v\v

\section{\bf  Moduli spaces of $J$-holomorphic maps}
\v
\subsection{Boundary conditions}

Consider the  symplectic manifold with cylindrical end
$$M^+=M^{+}_{0}\bigcup\left\{[0,\infty)\times \widetilde{M}\right\}.$$
Let $((\Sigma,{\bf j}); {\bf y}, {\bf p})$ be a connected semistable curve with $m$ marked points
${\bf y}=(y_1,...,y_m)$ and $\nu$ puncture points ${\bf p}=(p_1,...,p_{\nu})$,
and $u:{\Sigma} \rightarrow M^{+}$ be a
$J$-holomorphic map.  Let  $ \Sigma=\bigcup\limits_{v=1}^d(\Sigma_{v}, j_v) $ where $(\Sigma_{v},j_v)$ is a smooth Riemann surface and $\pi_{v}:\Sigma_{v} \rightarrow \Sigma$ is a continuous map.
 To describe the boundary conditions we consider two different cases separately:

\v
{\bf Case A .} Moduli space of $J$-holomorphic maps in contact geometry.
\v
Let $(\widetilde{M}, \lambda)$ be a compact contact manifold. Suppose that there exists a compact submanifold $\mathcal{F}\subset \widetilde{M}$  of dimension $\geq \;2$ satisfying
\v
{\bf (a)} $d\lambda|_{\mathcal{F}}=0,$
\v
{\bf (b)} every periodic orbit of the Reeb vector field lies in $\mathcal{F}$.
\v\n
Let $[c_{i}],i=1,\cdots,\fa$ be a bases of $H_{1}(\mathcal{F};\mathbb Z).$
\begin{definition}
Let ${\mathbf p}=(p_{1},\cdots,p_{\nu})$ be the order puncture points.
We assign a weight $\overrightarrow{\eta}$ to ${\mathbf p}$: \\
   $\overrightarrow{\eta}:{\mathbf p} \rightarrow   \mathbb Z_{>0}^{\oplus \fa} $ assigning a $\lambda_i=\sum_{l=1}^\fa \eta_{il}[c_{l}]$ to each puncture point $p_{i}$, where   $\eta_{il}\in \mathbb{Z}$.
   Choose the cylinder coordinates $(s_i,t_i)$ near $p_i$.
   We call a $J$-holomorphic map $u$ satisfies $(\overrightarrow{\eta})$ boundary condition if $u$ satisfies
 \begin{equation}
\lim\limits_{s_i\rightarrow \infty} u(s_i,S^1)\subset \mathcal{F},\;\;\forall\;1\leq i \leq \nu,
\end{equation}
\begin{equation}
[\lim\limits_{s_i\rightarrow \infty} u(s_i,S^1)]=\eta_i,\;\;\forall\;1\leq i \leq \nu.
\end{equation}
\end{definition}

\v
{\bf Case B.} Moduli space of $J$-holomorphic maps in $(M,L)$. \\
As we show in section \S\ref{neighbor_L} that $M-L$ can be considered as $M^+=M^{+}_{0}\bigcup\left\{[0,\infty)\times \widetilde{M}\right\}.$
Let $[c_{i}],i=1,\cdots,\fa$ is a bases in $H_{1}(L;\mathbb Z).$
\begin{definition} Let ${\mathbf p}=(p_{1},\cdots,p_{\nu})$ be the order puncture points.
We assign a weight $\overrightarrow{\mu}$ to ${\mathbf p}$: \\
   $\overrightarrow{\mu}:{\mathbf p} \rightarrow   \mathbb Z_{>0}^{\oplus \fa} $ assigning a $\mu_i=\sum_{l=1}^\fa \mu_{il}[c_{l}]$ to each puncture point $p_{i}$, where   $\mu_{il}\in \mathbb{Z}$.
   Choose the cylinder coordinates $(s_i,t_i)$ near $p_i$.
   We call a $J$-holomorphic map $u$ satisfies $(\overrightarrow{\mu})$ boundary condition if $u$ satisfies
 \begin{equation}
[\pi(\lim\limits_{s_i\rightarrow \infty} u(s_i,S^1))]=\mu_i,\;\;\forall\;1\leq i \leq \nu,
\end{equation}
where $\pi:T^{*}L\rightarrow L$ is  the canonical projection.
\end{definition}

\v

\v

\v

\subsection{Homology}

\v
{\bf Case A .} We fix $A\in H^{2}(M^{+},\mathcal{F};\mathbb Z)$ satisfying $\p A=\sum \eta_i.$ Consider a $J$-holomorphic map $u$ satisfying
\begin{equation}
[u(\Sigma)]=A.
\end{equation}
We show that the homology class $A$ give a bound of Energy. To simplify notation we assume that $\nu=1$. Let $(u, (\Sigma,{\bf j}),{\bf y},p)$ be a $J$-holomorphic map. By Theorem \ref{main theorem} $\tilde{u}$  converges to a $kT$-periodic orbit $x(kTt)$ as $z$ tends to $p$.  We construct a connected surface $W\subset \mathcal{F}$ with boundary $x(kTt)$.
Then $u(\Sigma)\cup W$ is a closed surface in $M^{+}$ and
$$[u(\Sigma)\cup W]\in  H^{2}(M^{+};\mathbb Z).$$
Denote $\bar{A}=[u(\Sigma)\cup W].$  By $d\lambda|_{\mathcal{F}}=0 $ we have
\begin{equation}\label{energy_bound}
 \omega({\bar{A}})=\int_{u(\Sigma)}\omega+\int_{W}\omega  =E_{\phi}(u)+\int_{W}d\lambda = E_{\phi}(u).
\end{equation}
 Let $W'\subset \mathcal{F}$ be another surface with boundary $x,$ denote $ \bar{A}'=[u(\Sigma)\cup W']\in  H^{2}(M^{+};\mathbb Z).$
We have $\omega(\bar{A})= \omega({\bar{A}}')=E_{\phi}(u),$ that is, $E_{\phi}(u)$ is independent of the choice of $W$ in $\mathcal{F}$.

\v
\v
{\bf Case B .} Let $A\in H^{2}(M,L;\mathbb Z)$ be a fixed homology class satisfying $\p A=\sum \mu_i.$ We have the same results.
\v\v

\section{\bf Compactness theorems }

\v

\subsection{Holomorphic blocks in $M^+$}\label{holomorphic_block_M}

Let $((\Sigma,{\bf j}); {\bf y}, {\bf p} )$ be a connected semistable curve with $m$ marked points
${\bf y}=(y_1,...,y_m)$ and $\nu $ puncture points ${\bf p} =(p_1 ,...,p_{\nu } )$.
Let  $u:{\Sigma} \rightarrow  M^{+}$ be a $J$-holomorphic map with finite energy. Suppose that $u(z)$ converges to a $k_i \cdot T_{i}$-periodic orbit $x(k_iT_{i}t)$ as $z$ tends to $p_i $.

\begin{definition}  A  J-holomorphic map $(u;((\Sigma,{\bf j}),{\bf y},{\bf p}))$ is said to be stable if for each $v$ one of the following conditions holds:
\begin{itemize}
\item[(1).] $u\circ \pi_{\Sigma_{v}}:\Sigma_{v}\rightarrow M^{+}$ is not a constant map.
\item[(2).] Let $val_v$ be the number of special points on $\Sigma_{v}$ which are nodal points, marked points or puncture points. Then
$val_v + 2g_{v}\geq 3.$
\end{itemize}
\end{definition}

\v

\begin{definition}\label{holomorphic_block_map_equiv_M}
Two   stable $J$-holomorphic maps $\Gamma =(u, (\Sigma,{\bf j}),{\bf y},{\bf p})$ and $\check{\Gamma} =(\check{u}, (\check{\Sigma},\check{\bf j}),{\bf \check{y}},{\bf \check{p}})$
is called equivalent if there exists a diffeomorphism $\varphi:\Sigma\rightarrow \check{\Sigma}$ such that  it can be lifted to bi-holomorphic isomorphisms $\varphi_{v w}:(\Sigma_{v},j_v)\rightarrow (\check\Sigma_{w},\check j_w)$ for
 each component $\Sigma_{v}$ of $\Sigma$, and
\begin{itemize}
\item[{\bf(1)}] $\varphi(y_i)= \check{y}_i$,  $\varphi(p_j)= \check{p_j}$ for any $1\leq i\leq m$, $1\leq j\leq \nu$,
\item[{\bf(2)}]  $\check{u}\circ \varphi= u$.
\end{itemize}
\end{definition}

\v

Denote by $\mathcal{M}_{A}(M,\mathcal{F};g,m+\nu,{\bf y},{\bf p},\overrightarrow{\mu})$ the moduli space of equivalence classes of all $J$-holomorphic curves in $M^{+}$ representing the homology class $A$ and satisfying $(\overrightarrow{\mu})$ boundary condition.

\begin{lemma}
There is a constant $C>0$ depending   on $A$ and $ \overrightarrow{\mu}$ such that for any $b=(u;(\Sigma,j),{\bf y},{\bf p})\in \mathcal{M}_{A}(M,\mathcal{F};g,m+\nu,{\bf y},{\bf p},\overrightarrow{\mu}) $ we have, over cylinder end,
\begin{equation}
E_{\phi}(u)+\int_{\Sigma}u^*d\lambda \leq C.
\end{equation}
\end{lemma}

The $D_{u}:\mathcal{W}^{1,p,\alpha}\rightarrow L^{p,\alpha}$ is a Fredholm operator with $ind=dim (ker D_{u})-dim(coker D_{u})$.
Put
$$Ind^L=ind +6(g-6)+2(m+\nu).$$
The virtual dimension of $\mathcal{M}_{A}(M,\mathcal{F};g,m+\nu,{\bf y},{\bf p},\overrightarrow{\mu})$ is $Ind^L$.
\vskip 0.1in
\noindent

\subsection{Holomorphic blocks in ${\mathbb{R}}\times \widetilde{M}$}\label{holomorphic_block_R}

In order to compactify the Moduli space $\mathcal{M}_{A}(M,\mathcal{F};g,m+\nu,{\bf y},{\bf p},\overrightarrow{\eta})$ and $\mathcal{M}_{A}(M,L;g,m+\nu,{\bf y},{\bf p},\overrightarrow{\mu})$ we need to consider $J$-holomorphic maps into
${\mathbb{R}}\times \widetilde{M}$.

Let $((\Sigma,{\bf j}); {\bf y}, {\bf p}^{+},{\bf p}^{-})$ be a connected semistable curve with $m$ marked points
${\bf y}=(y_1,...,y_m)$ and $\nu^{\pm}$ puncture points ${\bf p}^{+}=(p_1^{+},...,p_{\nu^{+}}^{+})$, ${\bf p}^{-}=(p_1^{-},...,p_{\nu^{-}}^{-})$,
and $u:{\Sigma} \rightarrow {\mathbb{R}}\times \widetilde{M}$ be a $J$-holomorphic map. Suppose that $u(z)$ converges to a $k_i^{\pm}\cdot T_{i^{\pm}}$-periodic orbit $x_{k_i^{\pm}}$ as $z$ tends to $p_i^{\pm}$.

\v

There is a $\mathbb R$ action, which induces a $\mathbb R$-action on the moduli space of $J$-holomorphic maps. We need mod this action.

\begin{definition}\label{equivalent_R}
Two $J$-holomorphic maps $\Gamma =(u, (\Sigma,{\bf j}),{\bf y}, {\bf p}^{+},{\bf p}^{-})$ and $\check{\Gamma} =(\check{u}, (\check{\Sigma},\check{\bf j}),{\bf \check{y}}, \check{\bf p}^{+},\check{\bf p}^{-})$
are called equivalent if there exists a diffeomorphism $\varphi:\Sigma\rightarrow \check{\Sigma}$ such that
  it can be lifted to bi-holomorphic isomorphisms $\varphi_{v \check{v}}:(\Sigma_{v},j_v)\rightarrow (\check\Sigma_{\check v},j_{\check v})$ for
 each component $\Sigma_{v}$ of $\Sigma$, and
\begin{itemize}
\item[{\bf(1)}] $\varphi(y_i)= \check{y_i}$,  $\varphi(p_j^{+})= \check{p_j}^{+}$, $\varphi(p_l^{-})= \check{p_l}^{-}$ for any $1\leq i\leq m$, $1\leq j\leq \nu^{+}$, $1\leq l\leq \nu^{-};$  $u$ and $\check u\circ \varphi$ converges to the same periodic orbit $x_{k_{i}^{\pm}}$ at $z$ tends to $p_{i}^{\pm};$
\item[{\bf(2)}] $\check{a}\circ \varphi= a+C,$ $\check{\tilde{u}}\circ \varphi=\tilde{u}$ for some constant $C$;
\end{itemize}
\end{definition}

\begin{definition}  A  J-holomorphic map $(u;((\Sigma,{\bf j}),{\bf y},{\bf p}))$ is said to be stable if for each $v$ one of the following conditions holds:
\begin{itemize}
\item[(1).] $\tilde{E}(u\circ \pi_{\Sigma_{v}})\ne 0$,
\item[(2).] Let $val_v$ be the number of special points on $\Sigma_{v}$ which are nodal points, marked points or puncture points. Then
$val_v + 2g_{v}\geq 3.$
\end{itemize}
\end{definition}

\v

For any $A\in H^{2}(\bR\times \widetilde M,\mathcal{F};\mathbb Z)$ we define $d\lambda (A)$ as following:
let $v: \mathbb R \times S^1\rightarrow \bR\times \widetilde M$ be a $C^{\infty}$ map such that $[v(\mathbb R \times S^1)]=A$, we define $d\lambda (A):= \int_{\mathbb R \times S^1} v^*(d\lambda).$

We fix $A\in H^{2}(\bR\times \widetilde M,\mathcal{F};\mathbb Z)$ and
${\bk}^\pm=(k_1^{\pm},...,k_{\nu^{\pm}}^{\pm})$
satisfying
\begin{equation}\label{energy_condtion_RM}
d\lambda(A)=\sum_{i=1}^{\nu^{+}} k_i^{+}\cdot T_{i^{+}}-\sum_{i=1}^{\nu^{-}} k_i^{-}\cdot T_{i^{-}} .
\end{equation}

We define
${\mathcal{M}}_{A}({\mathbb{R}}\times \widetilde{M},g,m+\nu^{+}+\nu^{-},{\bf k}^-, {\bf k}^+)$ to be the space of equivalence classes
of all stable $J$-holomorphic maps in ${\mathbb{R}}\times \widetilde{M}$ representing $A$ and converging to a $k_i^{\pm}\cdot T_{i^{\pm}}$-periodic orbits as $z$ tends to $p_i^{\pm}$. For any $(u,\Sigma,{\bf y},{\bf p}^{+},{\bf p}^{-})\in {\mathcal{M}}_{A}({\mathbb{R}}\times \widetilde{M},g,m+\nu^{+}+\nu^{-},{\bf k}^-, {\bf k}^+),$ by Stoke's formula we have
$$d\lambda(A)=\int_{\Sigma}u^*d\lambda=\sum_{i=1}^{\nu^{+}} \lambda(k_{i}^{+}x_{i}^{+}) -\sum_{i=1}^{\nu^{-}}\lambda(k_{i}^{-}x_{i}^{-}) =\sum_{i=1}^{\nu^{+}} k_i^{+}\cdot T_{i^{+}}-\sum_{i=1}^{\nu^{-}} k_i^{-}\cdot T_{i^{-}} .$$

We call ${\mathcal{M}}_{A}({\mathbb{R}}\times \widetilde{M};g,m+\nu^{+}+\nu^{-},{\bf k}^-, {\bf k}^+)$ a holomorphic rubber block in ${\mathbb{R}}\times \widetilde{M}$.
\v
Using the holomorphic blocks in $M^+$ and the holomorphic rubber blocks in ${\mathbb{R}}\times \widetilde{M}$ we can get the compactified moduli spaces $\overline{\mathcal{M}}_{A}(M,\mathcal{F};g,m+\nu,{\bf y},{\bf p},\overrightarrow{\eta})$ and $\overline{\mathcal{M}}_{A}(M,L;g,m+\nu,{\bf y},{\bf p},\overrightarrow{\mu})$ (see \cite{LS}).
\v\v

\section{\bf Contact invariants and Open string invariants}
\v

When the transversality fails we need to take the stabilization ( we use the terminology "regularization" in this paper). By a standard regularization procedure (see \cite{R2,CT,CLW-1,CLW-2,LR}) we get
a finite dimensional virtual  orbifold system
\[
\{(U_I,   E_I, \sigma_I)|  I\subset \{1, 2, \cdots, n\} \}
\]
  indexed by a   partially ordered
set $( I =2^{\{1, 2, \cdots, n\}}, \subset )$. Under some technical condition we can show that $\{\Sinverse\}$ is oriented, and the top strata of $U_{I}$ is a smooth orbifold.

Let $\Lambda=\{\Lambda_I\}$ be a partition of unity and $\{\Theta_I\}$ be a virtual  Euler form of $\{E_I\}$ such that
$\Lambda_I\Theta_I$ is compactly supported in $U_{I,\epsilon}$.
The contact invariant can be
defined as
\begin{equation}\label{c_invariant}
\Psi^{(C)}_{(A,g,m+\nu,\overrightarrow{\eta})}(\alpha_1,...,
\alpha_{m} ; \beta_{m + 1},..., \beta_{m+\nu})=\sum_I\int_{U_{I,\epsilon}}\prod_i e^*_i\alpha_i\wedge
\prod_j e^*_j\beta_j\wedge \Lambda_I \Theta_I.
\end{equation} for
$\alpha_i\in H^*(M^{+}, {\mathbb{R}})$ and $\beta_j\in H^*(\mathcal{F},
{\mathbb{R}})$ represented by differential form. Clearly, $\Psi^C=0$ if
$\sum \deg(\alpha_i)+ \sum \deg (\beta_i)\neq Ind^C$.
\v
Similarly, the open string invariant can be
defined as
\begin{equation}\label{l_invariant}
\Psi^{(L)}_{(A,g,m+\nu,\overrightarrow{\mu})}(\alpha_1,...,
\alpha_{m}; \beta_{m + 1},..., \beta_{m+\nu})=\sum_I\int_{U_{I,\epsilon}}\prod_i e^*_i\alpha_i\wedge
\prod_j e^*_j\beta_j\wedge \Lambda_I \Theta_I.
\end{equation} for
$\alpha_i\in H^*(M^{+}, {\mathbb{R}})$ and $\beta_j\in H^*(L,
{\mathbb{R}})$. Clearly, $\Psi^{(L)}=0$ if
$\sum \deg(\alpha_i)\neq Ind^L$.
\v
We consider the open manifold $\mathbb R\times \widetilde{M}$. By the same method above we can define the local contact invariants $\Psi^{(C,\ell)}_{(A,g,m+\nu,\overrightarrow{\eta})}$ and local open string invariants $\Psi^{(L,\ell)}_{(A,g,m+\nu,\overrightarrow{\mu})}$.

\v
It is proved that these integrals are independent of the choices of $\Theta_I$ and the regularization (see \cite{CT}), if they exist. The key issue is the convergence of the integrals near each lower strata.
\v
We use the gluing argument. We only consider the one nodal case, for general cases the proof is the same. Let $b=(u_1,u_2; \Sigma_1 \wedge \Sigma_2,j_1,j_2)$,
where $(\Sigma_1,j_1)$ and $(\Sigma_2,j_2)$ are smooth Riemann surfaces of genus $g_1$ and $g_2$ joining at $q$ and $u_i: \Sigma_i \rightarrow M$ are $J$-holomorphic with $u_1(q)=u_2(q).$
We use  the holomorphic cylindrical coordinates $(s_i, t_i)$ near $q$.  In terms of the holomorphic cylindrical
coordinates we write
$$\Sigma_1-\{q\}=\Sigma_{10}\bigcup\{[0,\infty)\times S^1\},$$
$$\Sigma_2-\{q\}=\Sigma_{20}\bigcup\{(-\infty,0]\times S^1\}.$$

\v

For any gluing parameter $r$ we construct a surface
$\Sigma_{r} =\Sigma_1 \#_{r} \Sigma_2 $. Then we glue the map $(u_1,u_2)$ to get the pregluing maps $u_{r}$, a family of approximate
$J$-holomorphic maps. Denote by $Q_b:L^{p,\alpha}\rightarrow K_{b}\times W^{1,p,\alpha}$ a right inverse of $D{\mathcal
S}_{b}.$  Then $D{\mathcal
S}_{b_{(r)}}$ is surjective for $r$ large enough. Moreover, there is a right inverses $Q_{b_{(r)}}$. By implicit function theorem, there exists a small open set $O$ of $0 \in Ker\;D{\mathcal S}_{b_{(r)}}$ and a unique smooth map
$$f_{(r)}: O\rightarrow
L^{p,\alpha}_{r}(u_{(r)}^{\ast}TM \otimes \wedge^{0,1})$$ such that for any $(\kappa,h)\in O$
$$ \bar{\p}_{J}\exp_{u_{(r)}}(\zeta)+\kappa_v=0,$$
where
 $$(\kappa_{v},\zeta)=((\kappa,h) + Q_{b_{(r)}}\circ f_{(r)}(\kappa,h)).$$
Put
$$E_1:=Ker D \mathcal{S}_{u_1},\;\;E_2:=Ker D \mathcal{S}_{u_2},\;\;\;\;\bh=  T_{q}M,$$
Denote
$$Ker D \mathcal{S}_{b}:=E_1\bigoplus_{\bh}E_2.$$
We can prove that when $r$ large enough there is an isomorphism
$$I_r: Ker D \mathcal{S}_{b}\longrightarrow Ker D \mathcal{S}_{b_{(r)}}.$$
Then we get a gluing map
$I_r(\kappa,\zeta)+Q_{b_{(r)}}\circ f_{(r)} \circ I_{r}(\kappa,\zeta)$ from $O_{j_1}\times O_{j_2} \times O$ into the moduli space, where $O$ is a neighborhood of $0$ in $Ker D \mathcal{S}_{b}$, $O_{j_k}$ is a neighborhood of $j_k$, $k=1,2$. Then $((j_1,j_2),r, \tau)$ is a local coordinate system in Delingne-Mumford space. We may choose $((j_1,j_2),r, \tau,\kappa, \zeta)$ as a local coordinate system in the moduli space. Using the same method in \cite{LS2} we can get the following

\begin{theorem}\label{coordinate_decay}
 Let $l\in \mathbb Z^+$ be a fixed integer. There exists positive  constants  $\mathcal C, \hbar, R_{0}$  and $ \alpha<\frac{1}{16e^{4}}$  such that for any $(\kappa,\zeta)\in Ker D \mathcal{S}_{b}$ with $\|(\kappa,\zeta)\|< \hbar$ the following estimates hold \\
 \begin{itemize}
\item[(I)] $\;\;\;\;\;
\left\|\frac{\partial }{\partial r}\left(I_r(\kappa,\zeta)+Q_{b_{(r)}}\circ f_{(r)} \circ I_{r}(\kappa,\zeta) \right) \right\|_{1,p,\alpha,r}\leq  \mathcal Ce^{-\alpha\tfrac{lr}{8} }.
$ \\

\item[(II)] $\;\;\;$ Restricting to the compact set $\{|s_i|\leq R_{0}\}$,  we have
$$
\left |\pi_{2}'\circ\frac{\partial }{\partial r}\left(I_r(\kappa,\zeta)+Q_{b_{(r)}}\circ f_{(r)} \circ I_{r}(\kappa,\zeta) \right) \right | \leq  \mathcal Ce^{-\alpha\tfrac{lr}{8} },
$$
 where  $\pi_{2}'(\kappa,h):=h$ for any $(\kappa,h)\in K_{b_{(r)}}\times \mathcal W_{r}^{1,p,\alpha}(\Sigma_{(r)};u_{(r)}
^{\ast}TM ).$
 \end{itemize}
  \end{theorem}
\vskip
0.1in \noindent
\begin{corollary}\label{exponential decay of map} Let $(\kappa,\zeta)$ be as in Theorem \ref{coordinate_decay}, denote
$$I_r(\kappa,\zeta)+Q_{b(r)}\circ f_{(r)}(I_r(\kappa,\zeta))=(\kappa_{r},h_{(r)}).$$
Then, restricting to the compact set $\{|s_i|\leq R_{0}\}$, we have
\begin{equation}
\left |\frac{\partial }{\partial r}\exp_{u_{(r)}}h_{(r)}\right |   \leq  \mathcal C e^{-\alpha\tfrac{lr}{8} }.
 \end{equation}
 %where $v:=\exp_{u_{(r)}}h_{(r)}$ and $P_{u_{(r)},v}$ denotes parallel transportation.
\end{corollary}
\v
We use this estimate to show the convergence of the integrals  \eqref{c_invariant}, \eqref{l_invariant}.

\v\v
\section{\bf Weinstein Conjecture}
\v 
The local contact invariants $\Psi^{(C,\ell)}_{(A,g,m+\nu,\overrightarrow{\eta})}$ and local open string invariants $\Psi^{(L,\ell)}_{(A,g,m+\nu,\overrightarrow{\mu})}$ can be used to study the Weinstein Conjecture. In particular, if $\Psi^{(C,\ell)}_{(A,g,m+\nu,\overrightarrow{\eta})}\ne 0$, there is a
perturbed $J$-holomorphic map $(\kappa,v)$ with finite energy, which satisfies
$$\bar{\partial}_{J}v + P_{u;v}\kappa =0.$$
Note that $\kappa$ supports in a compact set, so
there is a finite energy $J$-holomorphic map into $\mathbb R\times \widetilde{M}$. By the Theorem \ref{main theorem}, there is a periodic orbit of the Reeb vector field on $\widetilde{M}$. The argument holds also for $\Psi^{(L,\ell)}_{(A,g,m+\nu,\overrightarrow{\mu})}\ne 0$. Thus we have

\begin{theorem} The following hold
\begin{itemize}
\item[{\bf(1)}] If $\Psi^{(C,\ell)}_{(A,g,m+\nu,\overrightarrow{\eta})}\ne 0$ then there is a periodic orbit of the Reeb vector field on $\widetilde{M}$;
\item[{\bf(2)}] If $\Psi^{(L,\ell)}_{(A,g,m+\nu,\overrightarrow{\mu})}\ne 0$ then there is a periodic orbit of the Reeb vector field on $\widetilde{M}$.
\end{itemize}

\end{theorem}

\v

\end{document}